\documentclass[11pt]{asaproc}

\usepackage{graphicx}

\usepackage{amsthm, amstext, amsmath, amssymb}

\usepackage{txfonts}

\usepackage{float}

\usepackage{chngpage}
\usepackage{rotating}
\usepackage{graphicx}

\usepackage{subfig}

\usepackage{paralist}

\usepackage[%
    bookmarks=true,         
    unicode=false,          
    pdftoolbar=true,        
    pdfmenubar=true,        
    pdffitwindow=false,     
    pdfstartview={FitH},    
    pdftitle={Quickest Change-Point Detection: A Bird's Eye View},    
    pdfauthor={Aleksey S. Polunchenko, Grigory Sokolov and Wenyu Du},     
    pdfsubject={Quickest Change-Point Detection: A Bird's Eye View},   
    pdfcreator={Aleksey S. Polunchenko},   
    pdfproducer={Producer}, 
    pdfkeywords={CUSUM chart}{Quickest change detection}{Sequential analysis}{Sequential change-point detection}{Shiryaev's procedure}{Shiryaev--Roberts procedure}{Shiryaev--Roberts--Pollak procedure}{Shiryaev--Roberts--$r$ procedure}, 
    pdfnewwindow=true,      
    colorlinks=false,       
    linkcolor=red,          
    citecolor=green,        
    filecolor=magenta,      
    urlcolor=cyan           
]{hyperref}

\usepackage[T1]{fontenc}

\usepackage[authoryear, round]{natbib}

\renewcommand{\Pr}{\mathbb{P}} 
\DeclareMathOperator{\EV}{\mathbb{E}} 

\DeclareMathOperator{\LR}{\Lambda}

\DeclareMathOperator{\ADD}{ADD}

\DeclareMathOperator{\PFA}{PFA}
\DeclareMathOperator{\SADD}{SADD}

\DeclareMathOperator{\RIADD}{RIADD}
\DeclareMathOperator{\STADD}{STADD}
\DeclareMathOperator{\ESADD}{ESADD}

\DeclareMathOperator*{\esssup}{ess\,sup}

\newcommand{\T}{T}

\renewcommand{\le}{\leqslant} 
\renewcommand{\ge}{\geqslant}

\theoremstyle{plain} 
\newtheorem{theorem}{Theorem}

\newtheorem*{corollary}{Corollary}

\theoremstyle{definition} 
\newtheorem{definition}{Definition}

\graphicspath{{./gfx/}}

\usepackage{times}

\title{Quickest Change-Point Detection: A Bird's Eye View}

\author{Aleksey S. Polunchenko\thanks{Department of Mathematical Sciences, State University of New York at Binghamton, Binghamton, NY 13902--6000} \and Grigory Sokolov\thanks{Department of Mathematics, University of Southern California, Los Angeles, CA 90089--2532} \and Wenyu $\mathrm{Du}^{*}$} 

\begin{document}

\maketitle

\begin{abstract}
We provide a bird's eye view onto the area of sequential change-point detection. We focus on the discrete-time case with known pre- and post-change data distributions and offer a summary of the forefront asymptotic results established in each of the four major formulations of the underlying optimization problem: Bayesian, generalized Bayesian, minimax, and multi-cyclic.

\begin{keywords}
CUSUM chart,
Quickest change detection,
Sequential analysis,
Sequential change-point detection,
Shiryaev's procedure,
Shiryaev--Roberts procedure,
Shiryaev--Roberts--Pollak procedure,
Shiryaev--Roberts--$r$ procedure
\end{keywords}
\end{abstract}

\section{Introduction}
\label{sec:intro}

\vspace*{-90mm}
\begin{center}
{\small\it Proceedings of the 2013 Joint Statistical Meetings (JSM-2013)\\
Montr\'{e}al, Qu\'{e}bec, Canada, 3--8 August 2013}
\end{center}
\vspace*{76mm}
Quickest change-point detection is concerned with the design and analysis of procedures for ``on-the-go'' detection of possible changes in the characteristics of a running (random) process. Specifically, the process is assumed to be monitored continuously through sequentially made observations (e.g., measurements), and should their behavior suggest the process may have statistically changed, the aim is to conclude so within the fewest observations possible, subject to a tolerable level of the risk of false detection. See, e.g.,~\cite{Wald:Book47,Shiryaev:Book78,Siegmund:Book85,Poor+Hadjiliadis:Book08}. For nonparametric change-point detection theory see, e.g.,~\cite{Brodsky+Darkhovsky:Book93}. The area finds applications across many branches of science and engineering: industrial quality and process control (see, e.g.,~\citealp{Ryan:Book2011,Montgomery:Book2012,Wetherill+Brown:Book1991,Kenett+Zacks:Book1998,Shewhart:Book1931}), biostatistics (see, e.g.,~\citealp{Cohen:Book87}), clinical trials (see, e.g.,~\citealp{Siegmund:Book85}), econometrics (see, e.g.,~\citealp{Broemeling+Tsurumi:Book1987}), seismology (see, e.g.,~\citealp{Basseville+Nikiforov:Book93}), forensics, navigation, cybersecurity (see, e.g.,~\citealp{Tartakovsky+etal:SM2006-discussion} and \citealp{Polunchenko+etal:SA2012,Tartakovsy+etal:IEEE-JSTSP2013}), and communication systems (see, e.g.,~\citealp{Basseville+Nikiforov:Book93,Tartakovsky:Book91}) -- to name a few. See also, e.g.,~\cite{Chernoff:Book72}. A sequential change-point detection procedure, a rule whereby one stops and declares that (apparently) a change is in effect, is defined as a stopping time, $\T$, adapted to the observed data, $\{X_n\}_{n\ge1}$.

The desire to detect the change quickly causes one to be trigger-happy. That is, if one is too hasty, i.e., too quick to stop, the risk of a false detection is high. On the other hand, however, if one is too wary, i.e., too slow to stop, the delay to (correct) detection is substantial. Hence, there is a loss in either case and the essence of the problem is to attain a tradeoff between two contradicting performance measures -- the loss associated with the delay to detection of a true change and that associated with raising a false alarm. A good sequential detection policy is expected to minimize the average loss related to the detection delay, subject to a constraint on the loss associated with false alarms (or vice versa).

To put this idea on rigorous mathematical grounds one is to first formally define both the ``detection delay'' and the ``risk of raising a false alarm''. To this end, contemporary theory of sequential change-point detection distinguishes four different approaches: the minimax approach, the Bayesian approach, the generalized Bayesian approach, and the approach related to multi-cyclic detection of a distant change in a stationary regime. The aim of this paper is to give a brief overview of all four. For a more detailed overview see, e.g., \cite{Polunchenko+Tartakovsky:MCAP2012}, \cite{Tartakovsky+Moustakides:SA10}, and \cite{Tartakovsky+Veeravalli:TPA05}.

\section{Change-Point Models}
\label{sec:change-point-models} 

To formally state the general quickest change-point detection problem, one is to first introduce a change-point model, i.e., describe a probabilistic structure of the observations (independent, identically or non-identically distributed, correlated, etc.) as well as that of the change-point (unknown deterministic, random completely or partially dependent on the observed data, random fully independent from the observations). To this end, a myriad of scenarios is possible; see, e.g.,~\cite{Fuh:AS03,Fuh:AS04}, \cite{Tartakovsky:Book91,Tartakovsky:TPA09}, \cite{Tartakovsky+Moustakides:SA10},~\cite{Lai:JRSS95,Lai:IEEE-IT98}, \cite{Shiryaev:SMD61,Shiryaev:TPA63,Shiryaev:Book78,Shiryaev:TPA09,Shiryaev:SA10},~\cite{Tartakovsky+Veeravalli:TPA05}, and~\cite{Polunchenko+Tartakovsky:MCAP2012}. This section is intended to review the major ones.

Fix a probability triple $(\Omega,\mathcal{F},\Pr)$, where $\mathcal{F}=\vee_{n\ge 0}\mathcal{F}_n$, $\mathcal{F}_n$ is the sigma-algebra generated by the first $n\ge1$ observations ($\mathcal{F}_{0}=\{\varnothing,\Omega\}$ is the trivial sigma-algebra), and $\Pr\colon\mathcal{F}\mapsto[0,1]$ is a probability measure. Let $\Pr_{\infty}$ and $\Pr_{0}$ be two mutually locally absolutely continuous (i.e., equivalent) probability measures; for a general case with singular measures present see~\cite{Shiryaev:TPA09}. For $d=\{0,\infty\}$, write $\Pr_{d}^{(n)}=\Pr_{d}|_{\mathcal{F}_n}$ for the restriction of $\Pr_{d}$ to $\mathcal{F}_n$, and let $p_{d}^{(n)}(\cdot)$ be the density of $\Pr_{d}^{(n)}$ (with respect to a dominating sigma-finite measure).

Let $\{X_n\}_{n\ge1}$ be the series of observations such that $X_1,X_2,\ldots,X_{\nu}$, for some $\nu$, adhere to measure $\Pr_{\infty}$ (``normal'' regime), but $X_{\nu+1},X_{\nu+2},\ldots$ follow measure $\Pr_{0}$ (``abnormal'' regime). That is, at an unknown time instant $\nu$ (change-point), the observations undergo a change-of-regime from ``normal'' to ``abnormal''. Hence, $\nu$ is the serial number of the last normal observation, so that if $\nu=0$, then the entire series $\{X_n\}_{n\ge1}$ is in the abnormal regime admitting measure $\Pr_{0}$, while if $\nu=\infty$, then $\{X_n\}_{n\ge1}$ is in the normal regime admitting measure $\Pr_{\infty}$ (i.e., there is no change).

For every fixed $\nu\ge0$, the change-of-regime in the series $\{X_n\}_{n\ge1}$ generates a new probability measure $\Pr_{\nu}$. We will now construct the pdf $p_{\nu}^{(n)}(\boldsymbol{X}_1^n)$ of $\Pr_{\nu}^{(n)}$ for $n\ge1$ and $\nu\ge0$ in the most general case. For the sake of brevity, we will omit the superscript and write $p_{\nu}(\boldsymbol{X}_1^n)$.

For $1\le i\le j$, let $\boldsymbol{X}_i^j=(X_i,X_{i+1},\ldots,X_j)$, that is, $\boldsymbol{X}_i^j$ is a sample of $j-i+1$ successive observations indexed from $i$ through $j$. Hence, if the sample $\boldsymbol{X}_1^n=(X_1,X_{2},\ldots,X_n)$ is observed, then $\boldsymbol{X}_1^k=(X_1,\ldots,X_k)$ is the vector of the first $k$ observations in this sample and $\boldsymbol{X}_{k+1}^n=(X_{k+1},\ldots,X_n)$ is the vector of the rest of the observations in the sample, from $k+1$ to $n$.

First, suppose $\nu$ is deterministic unknown. This is the main assumption of the minimax approach. To get density $p_{\nu}(\boldsymbol{X}_1^n)$, observe that by the Bayes rule $p_{\infty}(\boldsymbol{X}_1^n)=p_{\infty}(\boldsymbol{X}_1^{\nu})\times p_{\infty}(\boldsymbol{X}_{\nu+1}^n|\boldsymbol{X}_1^{\nu})$ and $p_{0}(\boldsymbol{X}_1^n)=p_{0}(\boldsymbol{X}_1^{\nu})\times p_{0}(\boldsymbol{X}_{\nu+1}^n|\boldsymbol{X}_1^{\nu})$, whence by combining the first factor of the pre-change density, $p_{\infty}(\boldsymbol{X}_1^n)$, with the second one of the post-change density, $p_{0}(\boldsymbol{X}_1^n)$, we obtain $p_{\nu}(\boldsymbol{X}_1^n)=p_{\infty}(\boldsymbol{X}_1^{\nu})\times p_{0}(\boldsymbol{X}_{\nu+1}^n|\boldsymbol{X}_1^{\nu})$,
or, after some more algebra using the Bayes rule,
\begin{align}\label{eq:non-iid-model-def}
p_{\nu}(\boldsymbol{X}_1^n)
&=
\left(\,\prod_{j=1}^{\nu} p_{\infty}^{(j)}(X_j|\boldsymbol{X}_1^{j-1})\right)
\times
\left(\,\prod_{j=\nu+1}^n p_{0}^{(j)}(X_j|\boldsymbol{X}_1^{j-1})\right),
\end{align}
where $p_{\infty}^{(j)}(X_j|\boldsymbol{X}_1^{j-1})$ and $p_{0}^{(j)}(X_j|\boldsymbol{X}_1^{j-1})$ are the conditional densities of the $j$-th observation, $X_j$, given the past information $\boldsymbol{X}_1^{j-1}$, $j\ge1$. Note that in general these densities depend on $j$. Hereafter it is understood that $\prod_{j=k+1}^n p_{d}^{(j)}(X_j | \boldsymbol{X}_1^{j-1})=1$ for $k\ge n$.

Model~\eqref{eq:non-iid-model-def} is very general: it does not require the observations to be independent or homogeneous. Suppose now that $\{X_n\}_{n\ge1}$ are {\em independent} and such that $X_1,\ldots,X_{\nu}$ are each distributed according to a common density $f(x)$, while $X_{\nu+1},X_{\nu+2},\ldots$ each follow a common density $g(x)\not\equiv f(x)$. This is the simplest and most prevalent case. From now on it will be referred to as the {\em iid case}, or the {\em iid model}. In this case, model~\eqref{eq:non-iid-model-def} reduces to
\begin{align}\label{iidmodel}
p_{\nu}(\boldsymbol{X}_1^n)
&=
\left(\,\prod_{j=1}^{\nu} f(X_j)\right)
\times
\left(\,\prod_{j=\nu+1}^n g(X_j)\right),
\end{align}
and it will be referenced repeatedly throughout the paper.

If the change-point, $\nu$, is random, which is the ground assumption of the Bayesian approach, then any change-point model has to be supplied with a change-point's {\em prior distribution}. To this end, let $\pi_{0}=\Pr(\nu\le0)$ and $\pi_n=\Pr(\nu=n|\boldsymbol{X}_1^n)$, $n\ge1$, and observe that the series $\{\pi_n\}_{n\ge0}$ is $\{\mathcal{F}_n\}$-adapted. That is, the probability of the change occurring at time instance $\nu=k$ depends on $\boldsymbol{X}_1^k$, the observations' history accumulated up to (and including) time moment $k\ge1$. With the so defined prior distribution one can describe very general change-point models, including those that assume $\nu$ is a $\{\mathcal{F}_n\}$-adapted stopping time; see~\cite{Moustakides:AS08}.

To conclude this section, we note that when the probability series $\{\pi_n\}_{n\ge0}$ depends on the observed data $\{X_n\}_{n\ge1}$, it is argumentative whether $\{\pi_n\}_{n\ge0}$ can be referred to as the change-point's {\em prior} distribution: it can just as well be viewed as the change-point's {\em a posteriori} distribution. However, a deeper discussion of this subject is out of scope to this paper, and from now on, we will assume that $\{\pi_n\}_{n\ge0}$ do not depend on $\{X_n\}_{n\ge1}$, in which case it represents the ``true'' prior distribution.

\section{Overview of Optimality Criteria}
\label{sec:opt-criteria-overview} 

\subsection{Bayesian Formulation}

The signature assumption of the Bayesian formulation is that the change-point is a random variable with a prior distribution. This is instrumental in certain applications (see, e.g.,~\citealp{Shiryaev:METC06,Shiryaev:SA10} and \citealp{Tartakovsky+Veeravalli:TPA05}), but mostly of interest since the limiting versions of Bayesian solutions lead to optimal or asymptotically optimal procedures in more practical minimax problems.

Let $\{\pi_k\}_{k\ge0}$ be a prior distribution of the change-point, $\nu$, where $\pi_{0}=\Pr(\nu\le0)$ and $\pi_k=\Pr(\nu=k)$ for $k \ge 1$. From the Bayesian point of view, the risk of sounding a false alarm is reasonable to measure by the Probability of False Alarm (PFA), which is defined as
\begin{align}\label{PFAdef}
\PFA^{\pi}(\T)
&=
\Pr^\pi(\T\le\nu)=\sum_{k=1}^\infty \pi_k \Pr_k(T \le k),
\end{align}
where $\Pr^\pi(\mathcal{A}) = \sum\limits_{k=0}^\infty \pi_k \Pr_k(\mathcal{A})$ and the $\pi$ in the superscript emphasizes the dependence on the prior distribution.  Note that summation in~\eqref{PFAdef} is over $k \ge 1$ since by convention  $\Pr_k(\T\ge 1)=1$, so that $\Pr_k(\T\le 0)=0$. The most popular and practically reasonable way to benchmark the detection delay is through the Average Detection Delay (ADD), which is defined as
\begin{align}\label{ADDBayes}
\ADD^{\pi}(\T)
&=
\EV^\pi [\T-\nu|\T>\nu]=\EV^\pi[(\T-\nu)^+]/\Pr^\pi(\T>\nu),
\end{align}
where hereafter $x^+=\max\{0,x\}$ and $\EV^\pi$ denotes expectation with respect to $\Pr^\pi$.

We now formally define the notion of Bayesian optimality. Let $\Delta_\alpha=\{\T\colon\PFA^{\pi}(\T)\le\alpha\} $ be the class of detection procedures (stopping times) for which the PFA does not exceed a preset (desired) level $\alpha\in(0,1)$. Then under the Bayesian approach one's aim is to
\begin{align}\label{eq:Bayesian-opt-problem-def}
&\text{find $\T_{\mathrm{opt}}\in\Delta_{\alpha}$ such that $\ADD^{\pi}(\T_{\mathrm{opt}})=
\inf_{\T\in\Delta_{\alpha}}\ADD^{\pi}(\T)$ for every $\alpha\in(0,1)$}.
\end{align}

For the iid model~\eqref{iidmodel} and under the assumption that the change-point $\nu$ has a {\em geometric} prior distribution this problem was solved by~\cite{Shiryaev:SMD61,Shiryaev:TPA63,Shiryaev:Book78}. Specifically, Shiryaev assumed that $\nu$ is distributed according to the zero-modified geometric distribution
\begin{align}\label{eq:Shiryaev-geometric-prior-def}
\Pr(\nu<0)
&=
\pi
\;\;\text{and}\;\;
\Pr(\nu=n)=(1-\pi)p(1-p)^n,
\;\;n\ge0,
\end{align}
where $\pi\in[0,1)$ and $p\in(0,1)$. This is equivalent to choosing the series $\{\pi_n\}_{n\ge0}$ as $\pi_{0}=\Pr(\nu\le0)=\pi+(1-\pi)p$ and $\pi_n=\Pr(\nu=n)=(1-\pi)p(1-p)^n$, $n\ge1$.

Observe now that if $\alpha\ge1-\pi$, then problem~\eqref{eq:Bayesian-opt-problem-def} can be solved by simply stopping right away. This clearly is a trivial solution, since for this strategy the ADD is exactly zero, and $\PFA^{\pi}(\T)=\Pr(\nu>0)=1-\pi$, so that the constraint $\PFA^{\pi}(\T)\le\alpha$ is satisfied.  Therefore, assume that $\alpha<1-\pi$ and in this case, \cite{Shiryaev:SMD61,Shiryaev:TPA63,Shiryaev:Book78} proved that the optimal detection procedure is based on testing the posterior probability of the change currently being in effect, $\Pr(\nu<n|\mathcal{F}_n)$,  against a certain detection threshold. The procedure stops as soon as $\Pr(\nu<n|\mathcal{F}_n)$ exceed the threshold. This strategy is known as the Shiryaev procedure. To guarantee its strict optimality the detection threshold should be set so as to guarantee that the PFA is exactly equal to the selected level $\alpha$, which is rarely possible.

The Shiryaev procedure will play an important role in the sequel when considering non-Bayes criteria. It is more convenient to express Shiryaev's procedure through the average likelihood ratio (LR) statistic
\begin{align}\label{eq:Rnp-S-def}
R_{n,p}
&=
\frac{\pi}{(1-\pi)p}\prod_{j=1}^n\left(\frac{\LR_j}{1-p}\right)+\sum_{k=1}^n\prod_{j=k}^n \left(\frac{\LR_j}{1-p}\right),
\end{align}
where $\LR_n=g(X_n)/f(X_n)$ is the ``instantaneous'' LR for the $n$-th data point, $X_n$. Indeed, by using the Bayes rule, one can show that
\begin{align}\label{eq:posterior-pr-formula}
\Pr(\nu<n|\mathcal{F}_n)
&=
\frac{R_{n,p}}{R_{n,p}+1/p},
\end{align}
whence it is readily seen that ``thresholding'' the posterior probability $\Pr(\nu<n|\mathcal{F}_n)$ is the same as ``thresholding'' the process $\{R_{n,p}\}_{n\ge1}$.  Therefore, the Shiryaev detection procedure has the form
\begin{align}\label{eq:Shirst}
\T_{\mathrm{S}}(A)
&=
\inf\{n\ge1\colon R_{n,p}\ge A\},
\end{align}
and if $A=A_\alpha$ can be selected in such a way that the PFA is exactly equal to $\alpha$, i.e., $\PFA^{\pi}(\T_{\rm S}(A_\alpha))=\alpha$, then it is strictly optimal in the class $\Delta(\alpha)$, that is, $\inf_{\T\in\Delta(\alpha)}\ADD^\pi(\T)=\ADD^\pi(\T_{\mathrm{S}}(A_\alpha))$  for any $0<\alpha<1-\pi$. Note that Shiryaev's statistic $R_{n,p}$ can be rewritten in the recursive form
\begin{align}\label{eq:Rnp-S-recurrent-def}
R_{n,p}
&=
(1+R_{n-1,p})\frac{\LR_n}{1-p},
\;\; n\ge 1,
\;\;\text{with}\;\;
R_{0,p}=\frac{\pi}{(1-\pi)p}.
\end{align}

We also note that~\eqref{eq:Rnp-S-def} and~\eqref{eq:posterior-pr-formula} remain true under the geometric prior distribution~\eqref{eq:Shiryaev-geometric-prior-def} even in the general non-iid case~\eqref{eq:non-iid-model-def}, with $\LR_n= g(X_n|\boldsymbol{X}_1^{n-1})/f(X_n|\boldsymbol{X}_1^{n-1})$. However, in order for the recursion~\eqref{eq:Rnp-S-recurrent-def} to hold in this case, $\{\LR_n\}_{n\ge1}$ should be independent of the change-point.

As $p\to0$, where $p$ is the parameter of the geometric prior~\eqref{eq:Shiryaev-geometric-prior-def}, the Shiryaev detection statistic~\eqref{eq:Rnp-S-recurrent-def} converges to what is known as the {\em Shiryaev--Roberts (SR) detection statistic}. The latter is the basis for the so-called {\em SR procedure}. As we will see, the SR procedure is a ``bridge'' between all four different approaches to change-point detection mentioned above.

For a general asymptotic Bayesian change-point detection theory in discrete time see \cite{Tartakovsky+Veeravalli:TPA05}. Specifically, this work addresses the Bayesian approach assuming merely that the prior distribution is independent of the observations, and the overall conclusion is twofold:
\begin{inparaenum}[\itshape a\upshape)]
    \item the Shiryaev procedure is asymptotically (as $\alpha\to0$) optimal in a very broad class of change-point models and prior distributions, and
    \item depending on the behavior of the prior distribution at the right tail, the SR procedure may or may not be asymptotically optimal.
\end{inparaenum}
Specifically, if the tail is exponential, the SR procedure is not asymptotically optimal, though it is asymptotically optimal if the tail is heavy. When the prior distribution is arbitrary and depends on the observations, we are not aware of any strict or asymptotic optimality results.

\subsection{Generalized Bayesian Formulation}

The generalized Bayesian approach is the limiting case of the Bayesian formulation, presented in the preceding section. Specifically, in the generalized Bayesian approach the change-point $\nu$ is assumed to be a ``generalized'' random variable with a uniform (improper) prior distribution.

First, return to the Bayesian constrained minimization problem~\eqref{eq:Bayesian-opt-problem-def}. Specifically, consider the iid model~\eqref{iidmodel} and assume that the change-point $\nu$ is distributed according to zero-modified geometric distribution~\eqref{eq:Shiryaev-geometric-prior-def}. Then the Shiryaev  procedure defined in~\eqref{eq:Rnp-S-recurrent-def} and~\eqref{eq:Shirst} is optimal if the threshold $A=A_\alpha$ is chosen so that $\PFA^{\pi}(\T_{\rm S}(A_\alpha))=\alpha$. Suppose now that $\pi=0$ and $p\to0$; this is turning the geometric prior~\eqref{eq:Shiryaev-geometric-prior-def} to an improper uniform distribution. It can be seen that in this case $\{R_{n,p}\}_{n\ge0}$ becomes $\{R_{n,0}\}_{n\ge0}$, where $R_{0,0}=0$ and $R_{n,0}=(1+R_{n-1,0})\LR_n$, $n\ge1$ with $\LR_n=g(X_n)/f(X_n)$. The limit $\{R_{n,0}\}_{n\ge0}$ is known as the SR statistic, and is customarily denoted as $\{R_n\}_{n\ge0}$, i.e., $R_n=R_{n,0}$ for all $n\ge0$; in particular, note that $R_0=0$.

Next, when $\pi=0$ and $p\to0$ it can also be shown that
\begin{align}\label{eq:limits}
\frac{\Pr(\T>\nu)}{p}
&\rightarrow\EV_{\infty}[\T]
\;\;\text{and}\;\;
\frac{\EV[(\T-\nu)^+]}{p}\rightarrow\sum_{k=0}^\infty\EV_k[(\T-k)^+],
\end{align}
where $\T$ is an arbitrary stopping time. As a result, one may conjecture that the SR procedure minimizes the {\em Relative Integral Average Detection Delay} (RIADD)
\begin{align}\label{eq:RIADD-def}
\RIADD(\T)
&=
\frac{\sum_{k=0}^\infty\EV_k[(\T-k)^+]}{\EV_{\infty}[\T]}
\end{align}
over all detection procedures for which the {\em Average Run Length (ARL) to false alarm}, $\EV_{\infty}[\T]$, is no less than $\gamma>1$, an {\it a~priori} set level.

Let
\begin{align}\label{eq:class-Delta-ARL-def}
\Delta(\gamma)
&=
\bigl\{\T\colon\EV_{\infty}[\T]\ge\gamma\bigr\},
\end{align}
be the class of detection procedures (stopping times) for which the ARL to false alarm $\EV_{\infty}[\T]$ is ``no worse'' than $\gamma>1$. Then under the generalized Bayesian formulation one's goal is to
\begin{align}\label{eq:genBayesproblem}
&\text{find $\T_{\mathrm{opt}}\in\Delta(\gamma)$ such that $\RIADD(\T_{\mathrm{opt}})=
\inf_{\T\in\Delta(\gamma)}\RIADD(\T)$ for every $\gamma>1$}.
\end{align}

We have already hinted that this problem is solved by the SR procedure. This was formally demonstrated by~\cite{Pollak+Tartakovsky:SS09} in the discrete-time iid case, and by~\cite{Shiryaev:TPA63} and~\cite{Feiberg+Shiryaev:SD06} in continuous time for detecting a shift in the mean of a Brownian motion.

We conclude this subsection with two remarks. First, observe that if the assumption $\pi=0$ is replaced with $\pi=rp$, where $r\ge0$ is a fixed number, then, as $p\to0$, the Shiryaev statistic $\{R_{n,p}\}_{n\ge0}$ converges to $\{R_n^r\}_{n\ge0}$, where $R_n^r=(1+R_{n-1}^r)\LR_n$, $n\ge1$ with $R_0^r=r\ge0$. This is the so-called {\em Shiryaev--Roberts--$r$ (SR--$r$) detection statistic}, and it is the basis for the SR--$r$ detection procedure that starts from an arbitrary deterministic point $r$. This procedure is due to~\cite{Moustakides+etal:SS11}. The SR--$r$ procedure possesses certain minimax properties (cf.~\citealp{Polunchenko+Tartakovsky:AS10} and~\citealp{Tartakovsky+Polunchenko:IWAP10}). We will discuss this procedure at greater length later.

Secondly, though the generalized Bayesian formulation is the limiting (as $p\to0$) case of the Bayesian approach, it may also be equivalently re-interpreted as a completely different approach -- {\em multi-cyclic disorder detection in a stationary regime}. We will consider this approach in Subsection~\ref{sec:opt-criteria-overview:multi-cyclic}.

\subsection{Minimax formulation}
\label{sec:opt-criteria-overview:minimax}

Contrary to the Bayesian formulation the minimax approach posits that the change-point is an unknown not necessarily random number. Even if it is random its distribution is unknown. The minimax approach has multiple optimality criteria.

First minimax theory is due to~\cite{Lorden:AMS71} who proposed to measure the risk of raising a false alarm by the ARL to false alarm $\EV_{\infty}[\T]$. As far as the risk associated with detection delay is concerned, Lorden suggested to use the ``worst-worst-case'' ADD defined as
\begin{align*}
\ESADD(\T)
&=
\sup_{0 \le \nu<\infty}\biggl\{\esssup\EV_{\nu}[(\T-\nu)^+|\mathcal{F}_{\nu}]\biggr\}.
\end{align*}
 Lorden's minimax optimization problem seeks to
\begin{align}\label{eq:Lorden-minimax-problem}
&\text{find $\T_{\mathrm{opt}}\in\Delta(\gamma)$ such that $\ESADD(\T_{\mathrm{opt}})=
\inf_{\T\in\Delta(\gamma)}\ESADD(\T)$ for every $\gamma>1$} ,
\end{align}
where $\Delta(\gamma)$ is the class of detection procedures with the lower bound $\gamma$ on the ARL to false alarm defined in~\eqref{eq:class-Delta-ARL-def}.

For the iid scenario~\eqref{iidmodel},~\cite{Lorden:AMS71} showed that Page's~\citeyearpar{Page:B54} Cumulative Sum (CUSUM) procedure is first-order asymptotically minimax as $\gamma\to\infty$. For any $\gamma>1$, this problem was solved by~\cite{Moustakides:AS86}, who showed that CUSUM  is exactly optimal (see also~\cite{Ritov:AS90} who reestablished Moustakides'~\citeyearpar{Moustakides:AS86} finding using a different decision-theoretic argument).

Though the strict $\ESADD(\T)$-optimality of the CUSUM procedure is a strong result, it is more natural to construct a procedure that minimizes the average (conditional) detection delay, $\EV_{\nu}[\T-\nu|\T>\nu]$, for all $\nu\ge0$ simultaneously. As no such uniformly optimal procedure is possible,~\cite{Pollak:AS85} suggested to revise Lorden's version of minimax optimality by replacing $\ESADD(\T)$ with
\begin{align*}
\SADD(\T)
&=
\sup_{0\le\nu<\infty}\EV_{\nu}[\T-\nu|\T>\nu],
\end{align*}
the worst conditional expected detection delay. Thus, Pollak's version of the minimax optimization problem seeks to
\begin{align}\label{eq:Pollak-minimax-problem}
&\text{find $\T_{\mathrm{opt}}\in\Delta(\gamma)$ such that $\SADD(\T_{\mathrm{opt}})=
\inf_{\T\in\Delta(\gamma)}\SADD(\T)$ for every $\gamma>1$}.
\end{align}

It is our opinion that $\SADD(\T)$ is better suited for practical purposes for two reasons. First,  Lorden's criterion is effectively a double-minimax approach, and therefore, is overly pessimistic in the sense that $\SADD(\T)\le\ESADD(\T)$. Second, it is directly connected to the conventional decision theoretic approach --- the optimization problem \eqref{eq:Pollak-minimax-problem} can be solved by finding the least favorable prior distribution. More specifically, since by the general decision theory the minimax solution corresponds to the (generalized) Bayesian solution with the least favorable prior distribution, it can be shown that $\sup_{\pi}\ADD^{\pi}(\T)=\SADD(\T)$, where $\ADD^\pi(\T)$ is defined in \eqref{ADDBayes}. In addition, unlike Lorden's minimax problem~\eqref{eq:Lorden-minimax-problem}, Pollak's minimax problem~\eqref{eq:Pollak-minimax-problem} is still not solved.  For these reasons, from now on, when considering the minimax approach, we focus on Pollak's supremum ADD measure $\SADD(\T)$. Some light as to the possible solution (in the iid case) is shed in the work of~\cite{Polunchenko+Tartakovsky:AS10,Tartakovsky+Polunchenko:IWAP10}, and~\cite{Moustakides+etal:SS11}. A synopsis of the results is given in the sequel.

Yet another way to gauge the false alarm risk is through the worst local (conditional) probability of sounding a false alarm within a time ``window'' of a given length. As argued by~\cite{Tartakovsky:IEEE-CDC05,Tartakovsky:SA08-discussion}, in many surveillance applications (e.g., target detection) this may be a better option than the ARL to false alarm: the latter is more global. Specifically, the concern is that for a generic detection procedure, $\T$, the ARL to false alarm, $\EV_{\infty}[\T]$, is not an exhaustive measure of the false alarm risk, unless the $\Pr_{\infty}$-distribution of $\T$ is geometric (at least approximately); see~\cite{Tartakovsky:IEEE-CDC05,Tartakovsky:SA08-discussion}. The geometric distribution is characterized entirely by a single parameter, which a) uniquely determines $\EV_{\infty}[\T]$, and b) is uniquely determined by $\EV_{\infty}[\T]$. For the iid model~\eqref{iidmodel},~\cite{Tartakovsky+etal:IWAP08,Pollak+Tartakovsky:TVP09} showed that under mild assumptions the $\Pr_{\infty}$-distribution of the stopping times associated with detection schemes from a certain class is asymptotically (as $\gamma\to\infty$) exponential with parameter $1/\EV_{\infty}[\T]$; the convergence is in the $L^p$ sense, where $p\ge1$. The class includes all of the most popular procedures. Hence, for the iid model~\eqref{iidmodel}, the ARL to false alarm is an acceptable measure of the false alarm rate. However, for a general non-iid model this is not necessarily true. Hence, alternative measures of the false alarm rate are in order.
As a result, if $\T$ is geometric, one can evaluate $\Pr_{\infty}(k<\T\le k+m|\T>k)$ for  any $k\ge0$ (in fact, for all $k\ge0$ at once).
Specifically, let
\begin{align}\label{eq:class-Delta-sup-PFA-window-def}
\Delta_{\alpha}^m
&=
\biggl\{\T\colon\sup_{k\ge0}\Pr_{\infty}(k<\T\le k+m|\T>k)\le\alpha\biggr\},
\end{align}
be the class of detection procedures for which $\Pr_{\infty}(k<\T\le k+m|\T>k)$, the conditional probability of raising a false alarm inside a sliding window of $m\ge1$ observations  is ``no worse'' than a certain {\it a~priori} chosen level $\alpha\in(0,1)$. The size of the window $m$ may either be fixed or go to infinity as $\alpha\to0$.

As argued by~\cite{Tartakovsky:IEEE-CDC05}, in general, $\sup_{k}\Pr_{\infty}(k<\T\le k+m|\T>k)\le\alpha$ is a {\em stronger} condition than $\EV_{\infty}[\T]\ge\gamma$. Hence, in general, $\Delta_{\alpha}^m\subset\Delta(\gamma)$. See also~\cite{Tartakovsky:SA09-discussion}. For a specific example where the optimization problem~\eqref{eq:Pollak-minimax-problem} is solved in the class~\eqref{eq:class-Delta-sup-PFA-window-def} see~\cite{Polunchenko+Tartakovsky:MCAP2012}.

\subsection{Multi-cyclic detection of a disorder in a stationary regime}
\label{sec:opt-criteria-overview:multi-cyclic}

\begin{figure}
    \centering
    \subfloat[An example of the behavior of a process of interest as exhibited through the series of observations $\{X_n\}_{n\ge1}$.]{
        \includegraphics[width=0.9\textwidth]{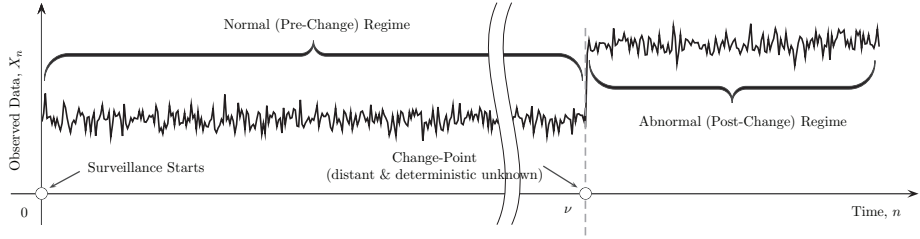}
    }\\ 
    \subfloat[An example of the behavior of the detection statistic when the decision to terminate surveillance is made {\em past} the change-point.]{
        \includegraphics[width=0.9\textwidth]{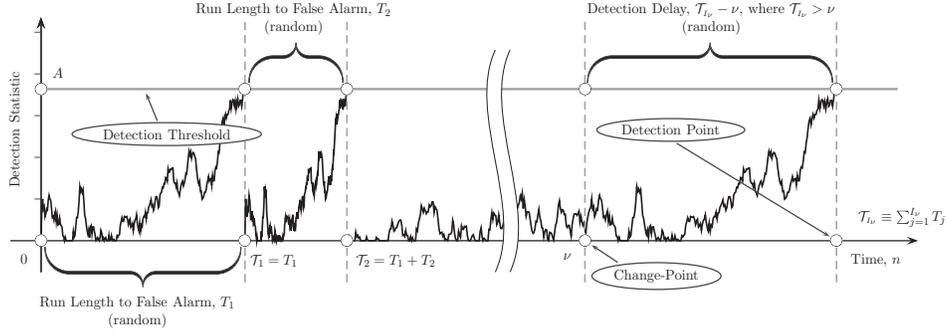}
    } 
    \caption{Multi-cyclic change-point detection in a stationary regime.}
    \label{fig:multi-cyclic-idea}
\end{figure}

Consider a context in which it is of utmost importance to detect the change as quickly as possible, even at the expense of raising many false alarms (using a repeated application of the same stopping rule) before the change occurs. This is equivalent to saying that the change-point $\nu$ is substantially larger than the tolerable level of false alarms $\gamma$. That is, the change ``strikes'' in a distant future and is preceded by a {\em stationary flow of false alarms}. This scenario is shown in~\autoref{fig:multi-cyclic-idea}. As one can see, the ARL to false alarm in this case is the mean time between (consecutive) false alarms, and therefore may be thought of the false alarm rate (or frequency).

As argued by~\cite{Pollak+Tartakovsky:SS09}, the multi-cyclic approach is instrumental in many surveillance applications, in particular in the areas concerned with intrusion/anomaly detection, e.g., cybersecurity and particularly detection of attacks in computer networks.

Formally, let $T_1,T_2,\ldots$ denote sequential independent repetitions of the same stopping time $\T$, and let $\mathcal{T}_{(j)}=T_{(1)}+T_{(2)}+\cdots+T_{(j)}$ be the time of the $j$-th alarm. Define $I_\nu=\min\{j\ge1\colon\mathcal{T}_{(j)}>\nu\}$. Put otherwise, $\mathcal{T}_{(I_\nu)}$ is the time of detection of the true change that occurs at the time instant $\nu$ after $I_\nu-1$ false alarms have been raised. Write
\begin{align*}
\STADD(\T)
&=
\lim_{\nu\to\infty}\EV_\nu[\mathcal{T}_{(I_\nu)}-\nu]
\end{align*}
for the limiting value of the ADD that we will refer to as the {\em stationary ADD} (STADD).

We now formally state the multi-cyclic change-point detection problem:
\begin{align}\label{eq:Multicycleoptim}
&\text{find $\T_{\mathrm{opt}}\in\Delta(\gamma)$ such that $\STADD(\T_{\mathrm{opt}})=
\inf_{\T\in\Delta(\gamma)}\STADD(\T)$ for every $\gamma>1$}
\end{align}
(among all multi-cyclic procedures).

For the iid model~\eqref{iidmodel}, this problem was solved by~\cite{Pollak+Tartakovsky:SS09}, who showed that the solution is the multi-cyclic SR procedure by arguing that $\STADD(\T)\equiv\RIADD(\T)$ defined in~\eqref{eq:RIADD-def}. This suggests that the optimal solution of the problem of multi-cyclic change-point detection in a stationary regime is completely equivalent to the solution of the generalized Bayesian problem. The exact result is stated in the next section.

\section{Optimality Properties of the Shiryaev--Roberts Detection Procedure}
\label{s:SRoptimality}

From now on we will confine ourselves to the iid scenario~\eqref{iidmodel}, i.e., assume that
\begin{inparaenum}[\itshape a\upshape)]
    \item the observations $\{X_n\}_{n\ge1}$ are independent throughout their history, and
    \item $X_1,\ldots,X_{\nu}$ are distributed according to a common known pdf $f(x)$ and $X_{\nu+1},X_{\nu+2},\ldots$ are distributed according to a common pdf $g(x)\not\equiv f(x)$, also known.
\end{inparaenum}

Let $\mathcal{H}_k\colon\nu=k$ for $0\le k<\infty$ and $\mathcal{H}_{\infty}\colon\nu=\infty$
be, respectively, the hypotheses that the change takes place at the time moment $\nu=k$, $k\ge0$, and that no change ever occurs. The densities of the sample $\boldsymbol{X}_1^n=(X_1,\ldots,X_n)$, $n\ge1$ under these hypotheses are given by
\begin{align*}
\begin{aligned}
p(\boldsymbol{X}_1^n|\mathcal{H}_{\infty})
&=
\prod_{j=1}^n f(X_j),
\;\;\text{and}\;\;
p(\boldsymbol{X}_1^n|\mathcal{H}_k)
=
\prod_{j=1}^k f(X_j)\prod_{j=k+1}^n
g(X_j)\;\;\text{for $k<n$},
\end{aligned}
\end{align*}
and $p(\boldsymbol{X}_1^n|\mathcal{H}_{\infty})=p(\boldsymbol{X}_1^n|\mathcal{H}_k)$ for $k\ge n$, so that the corresponding LR is
\begin{align*}
\LR_n^k
&=
\frac{p(\boldsymbol{X}_1^n|\mathcal{H}_k)}{p(\boldsymbol{X}_1^n|\mathcal{H}_{\infty})}
=\prod_{j=k+1}^n\LR_j\;\;\text{for}\;\; k<n,
\end{align*}
where $\LR_n=g(X_n)/f(X_n)$ is the ``instantaneous'' LR for the $n$-th observation $X_n$.

To decide in favor of one of the hypotheses $\mathcal{H}_k$ or $\mathcal{H}_{\infty}$, the likelihood ratios are ``fed'' to an appropriate sequential detection procedure, which is chosen according to the particular version of the optimization problem. In this section we are interested in the generalized Bayesian problem~\eqref{eq:genBayesproblem} and in the multi-cyclic disorder detection in a stationary regime~\eqref{eq:Multicycleoptim}. We have already remarked that for the iid model the SR procedure solves both of these problems. We preface the presentation of the exact results with the introduction of the SR procedure.

The SR procedure is due to the independent work of~\cite{Shiryaev:SMD61,Shiryaev:TPA63} and that of~\cite{Roberts:T66}. Specifically, Shiryaev considered the problem of detecting a change in the drift of a Brownian motion; Roberts focused on the case of detecting a shift in the mean of an iid Gaussian sequence. The name ``Shiryaev--Roberts'' was coined by~\cite{Pollak:AS85}. See~\cite{Pollak:IWSM09} for a brief account of the SR procedure's history.

Formally, the SR procedure is defined as the stopping time
\begin{align}\label{eq:T-SR-def}
\mathcal{S}_A
&=
\inf\{n\ge1\colon R_n\ge A\},
\end{align}
where $A>0$ is the detection threshold, and
\begin{align}\label{eq:Rn-SR-def}
R_n
&=
(1+R_{n-1})\LR_n,
\;\;n\ge1\;\;\text{with}\;\; R_0=0
\end{align}
is the SR detection statistic. As usual, we set  $\inf\{\varnothing\}=\infty$, i.e., $\mathcal{S}_A=\infty$ if $R_n$ never crosses $A$.

Recall first that $R_n=\lim_{p\to0}R_{n,p}$, where $R_{n,p}$ is the Shiryaev statistic given by recursion~\eqref{eq:Rnp-S-recurrent-def}. Recall also that the limiting relations~\eqref{eq:limits} hold. These allow us to conjecture that the SR procedure is optimal in the generalized Bayesian sense. In addition, since the RIADD is equal to the STADD of the multi-cyclic procedure, the repeated SR procedure should be optimal for detecting distant changes. The exact result is given next.

\begin{theorem}[\citealp{Pollak+Tartakovsky:SS09}]\label{Th:SRoptimality}
Let $\mathcal{S}_A$ be the SR procedure defined by \eqref{eq:T-SR-def} and~\eqref{eq:Rn-SR-def}. Suppose the detection threshold $A=A_{\gamma}$ is selected from the equation $\EV_{\infty}[\mathcal{S}_{A_{\gamma}}]=\gamma$, where $\gamma>1$ is the desired level of the ARL to false alarm.
\begin{enumerate}[\rm (i)]
    \item Then the SR procedure $\mathcal{S}_{A_{\gamma}}$ minimizes $\RIADD(\T)=\sum_{k=0}^\infty\EV_k[(\T-k)^+]/\EV_{\infty}[\T]$ over all stopping times $\T$ that satisfy $\EV_{\infty}[\T]\ge\gamma$, i.e., $\RIADD(\mathcal{S}_{A_{\gamma}})=\inf_{\T\in\Delta(\gamma)}\RIADD(\T)$ for every $\gamma>1$.
    \item Since $\RIADD(\T)\equiv\STADD(\T)$ for any stopping time $\T$, the SR procedure $\mathcal{S}_{A_{\gamma}}$ minimizes the stationary average detection delay among all multi-cyclic procedures in the class $\Delta(\gamma)$, i.e., $\STADD(\mathcal{S}_{A_{\gamma}})=\inf_{\T\in\Delta(\gamma)}\STADD(\T)$ for every $\gamma>1$.
\end{enumerate}
\end{theorem}

It is worth noting that the ARL to false alarm of the SR procedure satisfies the inequality $\EV_\infty[\mathcal{S}_{A}]\ge A$ for all $A>0$, which can be easily obtained by noticing that $R_n-n$ is a $\Pr_\infty$-martingale with mean zero. Also, asymptotically (as $A\to\infty$),  $\EV_\infty[\mathcal{S}_{A}]\approx  A/\zeta$, where the constant $0< \zeta <1$ is given by~\eqref{kappazeta} below (see~\citealp{Pollak:AS87}).  Hence, setting $A_\gamma=\gamma\zeta$ yields $\EV_\infty[\mathcal{S}_{A_\gamma}]\approx\gamma$, as $\gamma\to\infty$.

\section{Optimal and Nearly Optimal Minimax Detection Procedures}
\label{s:minimax-procedures} 

In this section, we will be concerned exclusively with the minimax problem in Pollak's setting~\eqref{eq:Pollak-minimax-problem}, assuming that the change-point $\nu$ is deterministic unknown. As of today, this problem is not solved in general. As has been indicated earlier, the usual way around this is to consider it asymptotically by allowing the ARL to false alarm $\gamma\to\infty$. The hope is to design such procedure $\T^*\in\Delta(\gamma)$ that $\SADD(\T^*)$ and the (unknown) optimum $\inf_{\T\in\Delta(\gamma)}\SADD(\T)$ will be in some sense ``close'' to each other in the limit, as $\gamma\to\infty$. To this end, the following three different types of asymptotic optimality are usually distinguished.

\begin{definition}[First-Order Asymptotic Optimality]
A procedure $\T^*\in\Delta(\gamma)$ is said to be {\em first-order asymptotically} optimal in the class $\Delta(\gamma)$ if $\SADD(\T^*)=\inf_{\T\in\Delta(\gamma)}\SADD(\T)[1+o(1)]$, as $\gamma\to\infty$, where from now on $o(1)\to0$, as $\gamma\to\infty$.
\end{definition}

\begin{definition}[Second-Order Asymptotic Optimality]
A procedure $\T^*\in\Delta(\gamma)$ is said to be {\em second-order asymptotically} optimal in the class $\Delta(\gamma)$  if $\SADD(\T^*)-\inf_{\T\in\Delta(\gamma)}\SADD(\T)=O(1)$, as $\gamma\to\infty$, where $O(1)$ stays bounded, as $\gamma\to\infty$.
\end{definition}

\begin{definition}[Third-Order Asymptotic Optimality]
A procedure $\T^*\in\Delta(\gamma)$ is said to be {\em third-order asymptotically} optimal in the class $\Delta(\gamma)$  if $\SADD(\T^*)-
\inf_{\T\in\Delta(\gamma)}\SADD(\T)=o(1)$, as $\gamma\to\infty$.
\end{definition}

\subsection{The Shiryaev--Roberts--Pollak procedure}

The question of what procedure minimizes Pollak's measure of detection delay $\SADD(\T)$ is an open issue. As an attempt to resolve the issue,~\cite{Pollak:AS85} proposed to ``tweak'' the SR procedure~\eqref{eq:T-SR-def}. This led to the new procedure that we will refer to as the Shiryaev--Roberts--Pollak (SRP) procedure. To facilitate the presentation of the latter, we first explain the heuristics.

As known from the general decision theory (see, e.g.,~\citealp[Theorem~2.11.3]{Ferguson:Book67}), an $\mathcal{F}_n$-adapted stopping time $\T$ solves~\eqref{eq:Pollak-minimax-problem} if
\begin{inparaenum}[\itshape a\upshape)]
     \item $\T$ is an extended Bayes rule,
     \item it is an equalizer, and
     \item it satisfies the false alarm constraint with equality.
\end{inparaenum}
A procedure is said to be an equalizer if its conditional risk (which we measure through $\EV_{\nu}[\T-\nu|\T>\nu]$) is constant for all $\nu\ge 0$, that is, $\EV_{0}[\T]=\EV_{\nu}[\T-\nu|\T>\nu]$ for all $\nu\ge1$. Of the three conditions the one that requires $\T$ to be an equalizer poses the most challenge.~\Citet{Pollak:AS85} came up with an elegant solution.

It turns out that the sequence $\EV_{\nu}[\mathcal{S}_A-\nu|\mathcal{S}_A>\nu]$ indexed by $\nu$ eventually {\em stabilizes}, i.e., it remains the same for all sufficiently large $\nu$. This happens because the SR detection statistic enters the quasi-stationary mode, which means that the conditional distribution $\Pr_{\infty}(R_n\le x|\mathcal{S}_A>n)$ no longer changes with time. If one could get to the quasi-stationary mode immediately, then the resulting procedure would have the same expected conditional detection delay for  all $\nu \ge 0$, i.e., it would be the equalizer. Thus, Pollak's~\citeyearpar{Pollak:AS85} idea was to start the SR detection statistic $\{R_n\}_{n\ge0}$, defined in~\eqref{eq:Rn-SR-def}, not from zero ($R_0=0$), but from a random point $R_0=R_0^Q$, where $R_0^Q$ is sampled from the {\em quasi-stationary distribution} of the SR statistic under the hypothesis $\mathcal{H}_{\infty}$ (which is a Markov Harris-recurrent process under $\mathcal{H}_{\infty}$). Specifically, the quasi-stationary cdf, $Q_A(x)$, is defined as
\begin{align}\label{eq:Q-SR-def}
Q_A(x)
&=
\lim_{n\to\infty}\Pr_{\infty}(R_n\le x|\mathcal{S}_A>n).
\end{align}

Therefore, the SRP procedure is defined as the stopping time
\begin{align}\label{eq:T-SRP-def}
\mathcal{S}_A^Q
&=
\inf\{n\ge1\colon R_n^Q\ge A\},
\end{align}
where $A>0$ is a detection threshold, and
\begin{align}\label{eq:R-SRP-def}
R_n^Q
&=
(1+R_{n-1}^Q)\LR_n, \;\; n\ge 1, \;\; R_0^Q\thicksim Q_A(x)
\end{align}
is the detection statistic.

We reiterate that, by design, the SRP procedure~\eqref{eq:T-SRP-def} and~\eqref{eq:R-SRP-def} is an equalizer: it delivers the same conditional average detection delay for any change-point $\nu\ge0$, that is, $\EV_{0}[\mathcal{S}_A^Q]=\EV_{\nu}[\mathcal{S}_A^Q-\nu|\mathcal{S}_A^Q>\nu]$ for all $\nu\ge1$. \Citet{Pollak:AS85} was able to demonstrate that the SRP procedure is third-order asymptotically optimal with respect to $\SADD(\T)$. We now state his result.
\begin{theorem}[\citealp{Pollak:AS85}]
Let $\EV_{0}[(\log\LR_1)^+]<\infty$. Suppose the detection threshold, $A$, of the SRP procedure, $\mathcal{S}_A^Q$, is set to the solution, $A_{\gamma}$, of the equation $\EV_{\infty}[\mathcal{S}_{A_\gamma}^Q]=\gamma$. Then $\SADD(\mathcal{S}_{A_\gamma}^Q)=\inf_{\T\in\Delta(\gamma)}\SADD(\T)+o(1)$, as $\gamma\to\infty$.
\end{theorem}
Recently,~\cite{Tartakovsky+etal:TPA2012} proved that $\EV_{0}[\mathcal{S}_{A}^Q]=(1/I)[\log A+\varkappa - C_{\infty}]+o(1)$,  as $A\to\infty$, provided $\EV_{0}[\log\LR_1]^2<\infty$,
where $\varkappa$ is the limiting average overshoot in the one-sided sequential test, which is a subject of renewal theory (see, e.g.,~\citealp{Woodroofe:Book82}), and $C_{\infty}$ is a constant that can be computed numerically (e.g., by Monte Carlo simulations). Both $\varkappa$ and $C_{\infty}$ are formally defined in the next subsection, where we reiterate the exact result of~\cite{Tartakovsky+etal:TPA2012}.

Note that for sufficiently large $\gamma$,
\begin{align}\label{ARLSRP}
\EV_{\infty}[\mathcal{S}_A^Q]
&\approx
(A/\zeta)-\mu_Q,
\;\;\text{where}\;\;
\mu_Q=\int_0^A y\,dQ_A(y),
\end{align}
i.e., $\mu_Q$ is the mean of the quasi-stationary distribution, and $\zeta$ is a constant defined in~\eqref{kappazeta} below. This approximation can be obtained by first noticing that for a fixed $R_0^Q=r$ the process $R_n^Q-r-n$ is a zero-mean $\Pr_\infty$-martingale, and then applying optional sampling theorem to this martingale as well as a renewal theoretic argument  (cf.~\citealp{Tartakovsky+etal:TPA2012}).

\subsection{The Shiryaev--Roberts--$r$ procedure}\label{ss:SRr}

Though the SRP procedure is practically appealing due to its third-order asymptotic optimality, it requires the knowledge of the quasi-stationary distribution~\eqref{eq:Q-SR-def} to implement. It is rare that this distribution can be expressed in a closed form; for examples where this is possible, see, e.g.,~\cite{Pollak:AS85},~\cite{Mevorach+Pollak:AJMMS91},~\cite{Polunchenko+Tartakovsky:AS10} and~\cite{Tartakovsky+Polunchenko:IWAP10}. As a result, the SRP procedure has not been used in practice.

To make the SRP procedure practical,~\cite{Moustakides+etal:SS11} proposed a numerical framework. More importantly,~\cite{Moustakides+etal:SS11} offered numerical evidence that there exist procedures that are uniformly better than the SRP procedure. Specifically, they regard starting off the original SR procedure at a fixed (but specially designed) $R_0=r$, $0\le r<A$, and defining the stopping time with this new deterministic initialization. Because of the importance of the starting point, they dubbed their procedure the SR--$r$ procedure.

Formally, the SR--$r$ procedure is defined as the stopping time
\begin{align}\label{eq:T-SR-r-def}
\mathcal{S}_A^r
&=
\inf\{n\ge1\colon R_n^r\ge A\},
\end{align}
where $A>0$ is the detection threshold, and
\begin{align}\label{eq:R-SR-r-def}
R_n^r
&=
(1+R_{n-1}^r)\LR_n,
\;\; n\ge 1,\;\;\text{with}\;\; R_0^r=r\ge0
\end{align}
is the SR--$r$ detection statistic.

\Citet{Moustakides+etal:SS11} show numerically that for certain values of the starting point, $R_0^r=r$, apparently, $\EV_{\nu}[\mathcal{S}_{A_1}^r-\nu|\mathcal{S}_{A_1}^r>\nu]$ is strictly less than $\EV_{\nu}[\mathcal{S}_{A_2}^Q-\nu|\mathcal{S}_{A_2}^Q>\nu]$ for all $\nu\ge 0$, where $A_1$ and $A_2$ are such that $\EV_{\infty}[\mathcal{S}_{A_1}^r]=\EV_{\infty}[\mathcal{S}_{A_2}^Q]$ (although the maximal expected delay is only slightly smaller for $\mathcal{S}_{A_1}^r$).

It turns out that using the ideas of~\cite{Moustakides+etal:SS11} we are able to design the initialization point $r=r(\gamma)$ in the SR--$r$ procedure~\eqref{eq:T-SR-r-def} so that this procedure is also third-order asymptotically optimal. In this respect, the average delay to detection at infinity $\ADD_{\infty}(\mathcal{S}_A^r)=\lim_{\nu\to\infty}\EV_{\nu}[\mathcal{S}_A^r-\nu|\mathcal{S}_A^r>\nu]$ plays the critical role. The following theorem, whose proof can be found in \cite{Polunchenko+Tartakovsky:AS10}, is important.

\begin{theorem}\label{Th1}
Let $\mathcal{S}_A^r$ be defined as in~\eqref{eq:T-SR-r-def} and~\eqref{eq:R-SR-r-def}, and let $A=A_\gamma$ be selected so that $\EV_{\infty}[\mathcal{S}_{A_\gamma}^r]=\gamma$. Then, for {\em every} $r\ge0$,
\begin{align}\label{IntADD}
\inf_{\T\in\Delta(\gamma)}\SADD(\T)
&\ge
\frac{r\EV_{0}[\mathcal{S}_{A_\gamma}^r]+\sum_{\nu=0}^\infty \EV_{\nu}[(\mathcal{S}_{A_\gamma}^r-\nu)^+]}{r+\EV_{\infty}[\mathcal{S}_{A_\gamma}^r]}=
\mathcal{J}_{\mathrm{B}}(\mathcal{S}_{A_\gamma}^r).
\end{align}
\end{theorem}

Note that~\autoref{Th1} suggests that if $r$ can be chosen so that the SR--$r$ procedure is an equalizer (i.e., $\EV_{0}[\mathcal{S}_A^r]=\EV_\nu[\mathcal{S}_A^r-\nu|\mathcal{S}_A^r>\nu]$ for all $\nu \ge 0$), then it is {\em exactly} optimal. This is because the right-hand side in~\eqref{IntADD} is equal to $\EV_{0}[\mathcal{S}_A^r]$, which, in turn, is equal to $\sup_\nu \EV_\nu[\mathcal{S}_A^r-\nu|\mathcal{S}_A^r>\nu]=\SADD(\mathcal{S}_A^r)$. Therefore, we have the following corollary.

\begin{corollary}\label{Cor:Cor1}
Let $A=A_\gamma$ be selected so that $\EV_{\infty}[\mathcal{S}_{A_\gamma}^r]=\gamma$. Assume that $r=r(\gamma)$ is chosen in such a way that the SR--$r$ procedure $\mathcal{S}_{A_\gamma}^{r(\gamma)}$ is an equalizer. Then it is strictly minimax in the class $\Delta(\gamma)$, i.e., $\inf_{\T\in\Delta(\gamma)}\SADD(\T)=
\SADD(\mathcal{S}_{A_\gamma}^{r(\gamma)})$.
\end{corollary}

\Citet{Polunchenko+Tartakovsky:AS10} and~\cite{Tartakovsky+Polunchenko:IWAP10} used this Corollary to prove that the SR--$r$ procedure with a specially designed $r=r_A$ is strictly optimal for two specific models. In general, \cite{Moustakides+etal:SS11} conjecture that the SR--$r$ procedure is third-order asymptotically minimax, and \cite{Tartakovsky+etal:TPA2012} show that this conjecture is true. We will state the exact result after we introduce some additional notation.

Let $S_n=\log\LR_1+\cdots+\log\LR_n$ and, for $a \ge 0$, introduce the one-sided stopping time $\tau_a=\inf\{n\ge1\colon S_n \ge a\}$. Let $\kappa_a=S_{\tau_a}-a$ be an overshoot (excess over the level $a$ at stopping), and let
\begin{align}\label{kappazeta}
\varkappa&=\lim_{a\to\infty}\EV_{0}[\kappa_a], \quad \zeta= \lim_{a\to\infty}\EV_{0}\left[e^{-\kappa_a}\right].
\end{align}

The constants $\varkappa>0$ and $0<\zeta<1$ depend on the model and can be computed numerically. Let $I=\EV_{0}[\log\LR_1]$ denote the Kullback--Leibler information number, and let $\tilde{V}_{\infty}=\sum_{j=1}^{\infty} e^{-S_j}$. Also, let $R_{\infty}$ be a random variable that has the $\Pr_{\infty}$-limiting (stationary) distribution of $R_n$, as $n \to\infty$, i.e., $Q_{\mathrm{ST}}(x)=\lim_{n\to\infty}\Pr_{\infty}(R_n\le x)=\Pr_{\infty}(R_{\infty}\le x)$. Let
\begin{align*}
C_{\infty}
&=
\EV[\log(1+R_{\infty}+\tilde{V}_{\infty})]
=
\int_0^\infty\int_0^\infty\log(1+x+y)\,dQ_{\mathrm{ST}}(x)\,d\tilde{Q}(y),
\end{align*}
where $\tilde{Q}(y)=\Pr_{0}(\tilde{V}_{\infty}\le y)$.

\begin{theorem}[\citealp{Tartakovsky+etal:TPA2012}]
Let $\EV_{0}[\log\LR_1]^2<\infty$ and let $\log\LR_1$ be non-arithmetic. Then the following assertions hold.
\begin{enumerate}[\rm (i)]
    \item $\inf_{\T\in\Delta(\gamma)}\SADD(\T)\ge(1/I)[\log(\gamma\zeta)+\varkappa-C_{\infty}]+o(1)$, as $\gamma \to \infty$.
    \item For any $r\ge0$,
\begin{align}\label{ADDinfty}
\ADD_{\infty}(\mathcal{S}_A^r)
&=
\EV_{0}[\mathcal{S}_A^{Q}]
=
\frac{1}{I}(\log A+\varkappa-C_{\infty})+o(1),
\;\;\text{as}\;\;A\to\infty.
\end{align}
\item Furthermore, if in the SR--$r$ procedure $A=A_\gamma=\gamma \zeta$ and the initialization point $r=o(\gamma)$ is selected so that $\SADD(\mathcal{S}_A^r)=\ADD_{\infty}(\mathcal{S}_A^r)$, then $\EV_{\infty}[\mathcal{S}_A^r] = \gamma (1+o(1))$ and $\SADD(\mathcal{S}_A^r)=(1/I)[\log(\gamma \zeta)+\varkappa-C_{\infty}]+o(1)$, as $\gamma\to\infty$.
\end{enumerate}
    Hence, the SR--$r$ procedure is third-order asymptotically optimal.
\end{theorem}

Also,
\begin{align}\label{SRrADD0}
\ADD_0(\mathcal{S}_A^r)
&=
\frac{1}{I}[\log A+\varkappa-C(r)]+o(1),
\;\;\text{as}\;\;
A\to\infty ,
\end{align}
where $C(r)=\EV[\log(1+ r+\tilde{V}_{\infty})]$. As we mentioned above, it is desirable to make the SR--$r$ procedure to look like equalizer by choosing the head start $r$, which can be achieved by equalizing $\ADD_0$ and $\ADD_\infty$. Comparing~\eqref{ADDinfty} and~\eqref{SRrADD0} we see that this property approximately holds  when $r$ is selected from the equation $C(r^*)=C_\infty$. This shows that asymptotically (as $\gamma\to\infty$) the ``optimal'' value $r^*$ is a fixed number that does not depend on $\gamma$. Clearly, this observation is important since it allows us to design the initialization point effectively and make the resulting procedure approximately optimal.

It is worth mentioning that $\SADD(\mathcal{S}_A)=\ADD_0(\mathcal{S}_A)=(1/I)[\log A+\varkappa-C(0)]+o(1)$, as $A\to\infty$, is true for the conventional SR procedure that starts from zero. Therefore, the SR procedure is only second-order asymptotically optimal. For sufficiently large $\gamma$, the difference between the supremum ADD-s of the SR procedure and the optimized SR--$r$ is given by $(C(0)-C_{\infty})/I$, which can be quite large if the Kullback--Leibler information number $I$ is small.

Note that similar to~\eqref{ARLSRP}, for sufficiently large $\gamma$, we have $\EV_{\infty}[\mathcal{S}_A^r]\approx(A/\zeta)- r$. For an example where distributions $Q_{\mathrm{ST}}(x)$ and $\tilde{Q}(x)$ and the constants $\varkappa$, $\zeta$, $C_{\infty}$, and $C(r)$ can be computed analytically see \cite{Polunchenko+Tartakovsky:MCAP2012}.

\bibliographystyle{apalike}
\bibliography{main,integral-equations,spc,operator-theory,size-bias}

\end{document}